\documentclass{article}
\usepackage{amssymb}
\usepackage{pdfsync}
\usepackage{color}

\newcommand{\dd}{{\rm \kern 3pt I\kern-9pt d}}

\newcommand{\Abar}{{\backslash\kern-8pt A}}

\newcommand{\CC}{\mbox{{\helv \i}{\rm\kern-5pt C}}}

\begin{document}
\begin{center}
\LARGE
The Lent Particle Method\\ for Marked Point Processes
\end{center}

\begin{flushright}
Nicolas Bouleau\\
Ecole des Ponts ParisTech
\end{flushright}
\vspace{.5cm}

Although introduced in the case of Poisson random measures (cf \cite{bouleau-denis}, \cite{bouleau-denis2}), the lent particle method applies as well in other situations. We study here the case of marked point processes. In this case the Malliavin calculus (here in the sense of Dirichlet forms) operates on the marks and the point process doesn't need to be Poisson. The proof of the method is even much simpler than in the case of Poisson random measures. We give applications to isotropic processes and to processes whose jumps are modified by independent diffusions.\\

\noindent{\bf 1.} {\bf Construction of the upper Dirichlet structure}\\

\noindent a) {\it Marked point processes}\\
Let $(X,\mathcal{X})$ and $(Y,\mathcal{Y})$ be two measurable spaces such that $\{x\}\in\mathcal{X}\;\forall x\in X$ and $\{y\}\in\mathcal{Y}\;\forall y\in Y$.

Let $\mathfrak{C}(X)$ be the configuration space of $X$ i.e. the space of countable sum $m$ of Dirac masses such that $m\{x\}\in\{0,1\}\;\forall x\in X$, so that $m$ may be indentified with its support. $\mathfrak{C}(X)$ is equipped with the smallest $\sigma$-field $\mathfrak{F}_X$ s.t. the maps $\omega\mapsto card(\omega\cap A)$ be measurable for any $A\in\mathcal{X}$.

Similarly we consider $\mathfrak{C}(X\times Y)$ equipped with $\mathfrak{F}_{X\times Y}$.

Let $\mu$ be a probability measure on $(Y,\mathcal{Y})$ and $\mathbb{Q}$ a probability measure on $(\mathfrak{C}(X),\mathfrak{F}_X)$. Let us denote by $M$ the random measure on $X$ with law  $\mathbb{Q}$.

For $F$ a function $\mathfrak{F}_{X\times Y}$-measurable and bounded, we may define a linear operator $S$ by putting
$$S(F)=\int F((x_1,y_1),\ldots,(x_n,y_n),\ldots)\,\mu(dy_1)\cdots\mu(dy_n)\cdots$$
the integral doesn't depend on the order of the numbering.
 $S(F)$ is $\mathfrak{F}_X$-measurable. Thus by
 $$\mathbb{P}(F)=\int S(F)\,d\mathbb{Q}$$
we define a probability measure on $(\mathfrak{C}(X\times Y),\mathfrak{F}_{X\times Y})$. We will say that $\mathbb{P}$ is the law of the random measure $M$ marked by $\mu$. It will be convenient to denote  $N=M\odot\mu$ this random measure of law $\mathbb{P}$.\\

\noindent b) {\it Dirichlet structure on a marked point process}\\
 We suppose that the measure $\mu$ is such that there exists a local Dirichlet structure with carr\'e du champ $(Y,\mathcal{Y},\mu,\mathbf{ d},\gamma)$. Although not necessary, we assume for simplicity
that constants belong to $\mathbf{ d}_{loc}$ (see Bouleau-Hirsch
\cite{bouleau-hirsch2} Chap. I Definition 7.1.3.) 
$$
1\in \mathbf{ d}_{loc} \makebox{ which implies }\ \gamma [1]=0.
$$
By the same argument as the theorem on products of Dirichlet structures (\cite{bouleau-hirsch2} Chap. V \S 2.2), the domain
$$\mathbb{D}=\{F\in L^2(\mathbb{P}),\textstyle \mbox{ for } \mathbb{Q}\mbox{-a.e.}\, m=\sum\varepsilon_{x_i}, \forall i, \mbox{ for }\mu\mbox{-a.e.}u_1,\ldots,\mu\mbox{-a.e.}u_{i-1},\mu\mbox{-a.e.}u_{i+1},\ldots $$ 
$$F((x_1,u_1),\ldots, (x_{i-1},u_{i-1}),(x_i,\,.\,),(x_{i+1},u_{i+1})\ldots)\in \mathbf{ d}$$
$$\mbox{ and }\mathbb{E}_\mathbb{P}[\sum_i(\gamma[F])(u_i)]<+\infty\}$$
and the operator $\Gamma[F]=\sum_i(\gamma[F])(u_i)$ define a local Dirichlet structure 

$$(\mathfrak{C}(X\times Y),\mathfrak{F}_{X\times Y}, \mathbb{P}, \mathbb{D}, \Gamma).$$

\noindent c) Let us recall the {\it Energy Image Density} property :  
For a $\sigma$-finite measure $\nu$ on some measurable space, a Dirichlet form on $L^2(\nu)$ with carr\'e du champ $\gamma$
is said to satisfy (EID) if for any $d$ and for any
$\mathbb{R}^d$-valued function $U$ whose components are in the
domain of the form
$$ U_*[({\det}\gamma[U,U^t])\cdot
\nu ]\ll \lambda^d $$ where $U_*$ denotes taking the image measure by $U$, {\rm det} denotes the determinant, and $\lambda^d$ the Lebesgue measure on $\mathbb{R}^d$.

For a local Dirichlet structure with carr\'e du champ, the above property is always
true for real-valued functions in the domain of the form (Bouleau \cite{bouleau1}, Bouleau-Hirsch
 \cite{bouleau-hirsch2} Chap. I \S7). It has been conjectured in 1986 (Bouleau-Hirsch
 \cite{bouleau-hirsch1} p251) that (EID) were true for any $\mathbb{R}^d$-valued function
 whose components are in the domain of the form for any local Dirichlet structure with carr\'e
  du champ. This has been shown for the Wiener space equipped with the Ornstein-Uhlenbeck form and
  for some other structures by Bouleau-Hirsch (cf. \cite{bouleau-hirsch2} Chap. II \S 5 and Chap. V example 2.2.4) and also for the Poisson space by A. Coquio \cite{coquio} when the intensity measure is the Lebesgue measure on an open set, and in more general cases in \cite{bouleau-denis} thanks to a result of Song \cite{song}. But this conjecture
  being at present neither refuted nor proved in full generality, the property has to be established in every
  particular setting.\\

\noindent{\bf Lemma 1}. {\it If the structure $(Y,\mathcal{Y},\mu,\mathbf{ d},\gamma)$ is such that any finite product $(Y,\mathcal{Y},\mu,\mathbf{ d},\gamma)^n$, $n\in\mathbb{N}$, satisfies} (EID) {\it then the structure $(\mathfrak{C}(X\times Y),\mathfrak{F}_{X\times Y}, \mathbb{P}, \mathbb{D}, \Gamma)$ satisfies} (EID).\\

\noindent This is an application of prop. 2.2.3 and thm 2.2.1 of Chap. V of \cite{bouleau-hirsch2}.\\

\noindent d) {\it The lent particle method}.

Let us denote $\varpi$ the current point of the space  $\mathfrak{C}(X\times Y)$, and let us introduce the operators
$$\varepsilon^+_{(x,u)}\varpi=\varpi\cup\{(x,u)\}\qquad\varepsilon^-_{(x,u)}\varpi=\varpi\cap\{(x,u)\}^c$$
then we have the lent particle formula
\begin{equation}\forall F\in \mathbb{D}\qquad\Gamma[F]=\int\varepsilon^-\gamma\varepsilon^+F\;dN\end{equation}
\noindent{\it Proof}. For $F\in \mathbb{D}$ we have
$$
\begin{array}{rl}
\varepsilon^+F=&F((x,u),(x_1,u_1),\ldots, (x_i,u_i),\ldots)\\
\gamma\varepsilon^+F=&\gamma[F((x,.),(x_1,u_1),\ldots, (x_i,u_i),\ldots)](u)
\end{array}$$
and $\int\varepsilon^-\gamma\varepsilon^+F\;dN$ is the sum, when $(x,u)$ varies among the points $(x_i,u_i)\in\varpi$ of the preceding result. This makes
$$\sum_i\gamma_i[F],$$ exactly what we obtained by the product construction. This shows also, by the definition of $\mathbb{D}$, that the integral $\int\varepsilon^-\gamma\varepsilon^+F\;dN$ exists and belongs to $L^1(\mathbb{P})$. \hfill$\square$\\

\noindent e) {\it Gradient}. Let us explain how could be done the construction of a gradient for the structure $(\mathfrak{C}(X\times Y),\mathfrak{F}_{X\times Y}, \mathbb{P}, \mathbb{D}, \Gamma)$ starting from a gradient for the structure $(Y,\mathcal{Y},\mu,\mathbf{ d},\gamma)$.

Let us suppose that the structure $(Y,\mathcal{Y},\mu,\mathbf{ d},\gamma)$ is such that the Hilbert space $\mathbf{d}$ be separable.
Then by a result of Mokobodzki (see Bouleau-Hirsch \cite{bouleau-hirsch2}, ex.5.9 p. 242) this
Dirichlet structure admits a gradient operator in the sense that
there exists a separable Hilbert space $H$ and a continuous linear
map $D$ from $\mathbf{ d}$ into $L^2 (Y,\mu;H)$ such that
\begin{itemize}
\item $\forall u\in \mathbf{ d}$, $\| D[u ]\|^2_H =\gamma[u]$. \item If
$F:\mathbb{R}\rightarrow \mathbb{R}$ is Lipschitz  then
\[\forall u\in \mathbf{ d},\ D[F\circ u]=(F'\circ u )Du.\]
\item If $F$ is $\mathcal{C}^1$ (continuously differentiable) and Lipschitz from $\mathbb{R}^d$ into
$\mathbb{R}$ (with $d\in \mathbb{N}$) then
\[ \forall u=(u_1 ,\cdots ,u_d) \in \mathbf{ d}^d ,\ D[F\circ
u]=\sum_{i=1}^d (F'_i \circ u ) D[u_i ].\]
\end{itemize}
As only the Hilbertian structure of $H$ plays a role, we can choose for
$H$
 a  space $L^2 (R,\mathcal{R} ,\rho)$ where $(R,\mathcal{R} ,\rho)$ is a
 probability space such that the dimension of the vector space $L^2 (R,\mathcal{R}
 ,\rho)$ is infinite. As usual, we identify $L^2 (\mu ;H)$ and
 $L^2 (Y\times R ,\mathcal{Y} \otimes \mathcal{R} ,\mu\times\rho)$ and we denote
 the gradient $D$ by  $\flat$:
 \[\forall u\in \mathbf{ d},\ Du=u^{\flat} \in L^2 (Y\times R ,\mathcal{Y} \otimes \mathcal{R}
 ,\mu\times\rho).\]
 Without loss of generality, we assume moreover that operator
 $\flat$ takes its values in the orthogonal space of $1$ in $L^2
 (R,\mathcal{R} ,\rho)$, in other words we take for $H$ the orthogonal of
 $1$. So that we have
$$
 \forall u\in \mathbf{ d} ,\ \int u^{\flat}d\rho =0\ \ \mu\mbox{-}a.e.
$$
Finally,  by the hypothesis on $\gamma$ we have
$$
1\in \mathbf{ d}_{loc} \makebox{ which implies }\ \gamma [1]=0 \makebox{
and  } 1^{\flat}=0.
$$
With these tools and hypotheses we obtain easily a gradient for the structure  $(\mathfrak{C}(X\times Y),\mathfrak{F}_{X\times Y}, \mathbb{P}, \mathbb{D}, \Gamma)$. We have to follow the same construction as above replacing the measure $\mathbb{Q}\times\mu^{\mathbb{N}}$ by the measure $\mathbb{Q}\times\mu^{\mathbb{N}}\times\rho^{\mathbb{N}}$. This yields a random measure $N\odot\rho=M\odot\mu\times\rho$ defined under the probability measure $\mathbb{P}\times\rho^\mathbb{N}$.

Now it is straightforward to show that the formula 
$$F^\sharp=\int \varepsilon^-(\varepsilon^+F)^\flat\;dN\odot\rho$$ for $F\in\mathbb{D}$ defines a gradient for the structure $(\mathfrak{C}(X\times Y),\mathfrak{F}_{X\times Y}, \mathbb{P}, \mathbb{D}, \Gamma)$ with values in $L^2(\mathbb{P}\times\rho^\mathbb{N})$. The existence of the integral $\int\varepsilon^-(\varepsilon^+F)^\flat\;dN\odot\rho$ comes from the fact that it is controlled by that of $\int\varepsilon^-\gamma\varepsilon^+F\;dN$ thanks to
$$\rho^\mathbb{N}\left\{(\int\varepsilon^-(\varepsilon^+F)^\flat\;dN\odot\rho)^2\right\}=
\int\int(\varepsilon^-(\varepsilon^+F)^\flat)^2d\rho dN=\int\varepsilon^-\gamma[\varepsilon^+F]dN$$ (similar formula as in Corollary 12 of \cite{bouleau-denis}).

\noindent{\it Example. } If $F=e^{-N(f)}$, then $$\varepsilon^+_{(x,u)}F=e^{-N(f)}e^{-f(x,u)}$$  
$$\gamma\varepsilon^+_{(x,u)}F=e^{-2N(f)}e^{-2f(x,u)}\gamma[f]$$
$$\int\varepsilon^-\gamma\varepsilon^+F\;dN=e^{-2N(f)}N(\gamma[f])\qquad(=e^{-2N(f)}\Gamma[N(f)])$$

( $\Gamma[N(f)]=N(\gamma[f])$ even in the non Poissonian case).\\

Let us summarize this construction which gives a result, similar to Theorem 17 of \cite{bouleau-denis}, obtained much more easily here for marked point processes than for random Poisson measures.

\noindent{\bf Theorem 2.} {\it The carr\'e du champ operator of the upper Dirichlet structure $(\mathfrak{C}(X\times Y),\mathfrak{F}_{X\times Y}, \mathbb{P}, \mathbb{D}, \Gamma)$ satisfies $\forall F\in\mathbb{D}$
$$\Gamma[F]=\int\varepsilon^-\gamma[\varepsilon^+ F]dN$$ and this structure satisfies {\rm (EID)} as soon as every finite product $(Y,\mathcal{Y},\mu,\mathbf{ d},\gamma)^n$ satisfies {\rm (EID)}.}\\

\noindent{\bf 2.  Application to isotropic processes.}

 Let us consider a L\'evy process $Z=(Z^1,Z^2)$ with values in $\mathbb{R}^2$ and L\'evy measure $\sigma(dx,dy)=\nu(dr)\tau(d\theta)$ where $\tau$ is the uniform probability on the circle. Let us suppose that  $Z$ is centered without Gaussian part and that $\sigma$ integrates $r^2=x^2+y^2$. Let $N$ be the Poisson measure such that for any $h_1$ and $h_2$ in $L^2(ds)$
$$\int_0^t\;h_1(s)dZ^1_s+h_2(s)dZ^2_s=\int1_{[0,t]}(s)(h_1(s)x+h_2(s)y)\,\tilde{N}(dsdxdy).$$
Let us construct the upper Dirichlet structure starting from the classical structure on the unit circle with domain $H^1$. And let us consider as illustration the very simple functional
 $F=Z_t=(r_t\cos\theta_t, r_t\sin\theta_t)$
$$\varepsilon^+_{(t_0,r_0,\theta_0)}F=(Z^1_t+1_{t\geq t_0}r_0\cos\theta_0,Z^2_t+1_{t\geq t_0}r_0\sin\theta_0)$$
$$
\gamma\,\varepsilon^+F=1_{t\geq t_0}\left(
\begin{array}{lr}
\sin^2\theta_0&\cos\theta_0\sin\theta_0\\
\cos\theta_0\sin\theta_0&\cos^2\theta_0
\end{array}
\right)r^2_0
$$
$$\Gamma[F]=\int\varepsilon^-\gamma\,\varepsilon^+F\,dN=\int_0^tr^2\left(
\begin{array}{lr}
\sin^2\theta&\cos\theta\sin\theta\\
\cos\theta\sin\theta&\cos^2\theta
\end{array}\right)
\,N(dsdrd\theta).$$
As soon as $\nu$ has an infinite mass, $\forall t>0, \exists r_1\neq 0, r_2\neq 0$ et $\theta_1\neq\theta_2$ s.t. 
$$\Gamma[F]\geq r^2_1\wedge r^2_2
\left(
\begin{array}{lr}
\sin^2\theta_1+\sin^2\theta_2&\cos\theta_1\sin\theta_1+\cos\theta_2\sin\theta_2\\
\cos\theta_1\sin\theta_1+\cos\theta_2\sin\theta_2&\cos^2\theta_1+\cos^2\theta_2
\end{array}
\right)
$$ in the sense of positive symmetric matrices. Hence it follows that
$$\det\Gamma[F]\geq (r^2_1\wedge r^2_2)^2\sin^2(\theta_1-\theta_2)>0.$$
So that $Z_t$ possesses a density on $\mathbb{R}^2$, as soon as $\nu(\mathbb{R}_+^\ast)=+\infty$. This result is probably known although not contained in the criterion of Sato \cite{sato} which supposes $\nu$ absolutely continuous. (Here $\nu$ may be possibly a weighted sum of Dirac masses because it doesn't carry any Dirichlet form).

The measure on the circle need not to be uniform provided that it carries a Dirichlet form such that its $n$-th powers satisfy (EID).
The idea generalizes obviously replacing the circle by a $d$-dimensional sphere.

Actually, the process $Z$ doesn't need to be L\'evy. The method applies as well for instance to  a real process purely discountinuous  if we modify its jumps by i.i.d. transformations.\\

\noindent{\bf 3. Insight on transform of L\'evy processes by diffusions.} 

Since the Wiener measure is a probability measure we may take for  $(Y,\mathcal{Y},\mu)$ the Wiener space equipped with the Ornstein-Uhlenbeck structure. We know that  (EID) is fulfilled as asked in Thm 2.

Let us consider the SDE 
\begin{equation}\label{eds}X_t^x=x+\sum_{j=1}^d\int_0^t A_j(X_\tau^x,x) dB_\tau^j+\int_0^tB(X_\tau^x,x)d\tau\end{equation} 
where $x\in\mathbb{R}^m$. The coefficients are $\mathcal{C}^1\cap Lip$ with respect to the first argument.

Let us take for $(X,\mathcal{X})$  the Euclidean space $(\mathbb{R}_+\!\times\mathbb{R}^m,\mathcal{B}(\mathbb{R}_+\!\times\mathbb{R}^m))$.
Let $M$ be a random Poisson measure on $\mathbb{R}_+\times\mathbb{R}^m$ with intensity $ds\times\nu$ and law  $\mathbb{Q}$ associated with a L\'evy process $Z$. We denote $\varpi=\sum_\alpha \varepsilon_{(s_\alpha,x_\alpha)}$ the current point of $\mathfrak{C}(X)$.\\

Equation (\ref{eds}) is not that of a homogeneous Markov process because of the second argument in the coefficients. We can nevertheless define $\Pi_{t,x}(d\xi)$ to be the law of $X^x_t$ and $\nu\Pi_t=\int\nu(dx)\Pi_{t,x}$ to be the law of $X_t$ starting with the measure $\nu$.\\

\noindent{\bf Lemma 3.} {\it If the coefficients $A_j$, $B$ are Lipschitz with respect to the first argument with constant independent of $x$ and vanish at zero, the transition  $\Pi_t$  preserves L\'evy measures and measures integrating $x\mapsto |x|\wedge 1$.}

\noindent{\it Proof.} By Gronwall lemma for p=1 or p=2, $\mathbb{E}|X_t^x|^p\leq k|x|^pe^{kt}$, this means that $\nu \Pi_t$ is a L\'evy measure for any L\'evy measure $\nu$ and the lemma. \hfill$\Box$\\

The  transformed L\'evy process $(T_t(Z))_s$  whose jumps are modified independently by the diffusion (\ref{eds}), which is a L\'evy process with L\'evy measure $\nu\Pi_t$, is a functional $F$ of the marked point process. Let us suppose for simplicity that the jumps of $Z$ are summable, i.e. that $\nu$ integrates $x\mapsto|x|\wedge 1$, then $F$ may be written
$$ F=\int_{[0,s]\times\mathbb{R}^m\times Y}X_t^x(y)N(dsdxdy)$$ with as above $N=M\odot\mu$.
The lent particle formula gives
$$F^\sharp=\int_{[0,s]\times\mathbb{R}^m\times Y\times R}(X_t^x)^\flat \;d(N\odot\rho)$$ 
and
$$\Gamma[F]=\int_{[0,s]\times\mathbb{R}^m\times Y}\gamma[X_t^x] \;dN.$$
 Now $ (X_t^x)^\flat$ and $\gamma[X_t^x]$ are known by the usual Malliavin calculus  : $(.)^\flat$ is a gradient on the Wiener space associated with the O-U structure, for which we can choose (cf \cite{bouleau-hirsch2}) the operator  defined by
$$(\int h(s)dB_s^j)^\flat=\int h(s)d\hat{B}_s^j\quad h\in L^2(\mathbb{R}_+)$$ where $\hat{B}^j$ are independent copies of $B^j$. 

$$(X_t^x)^\flat=K_t\int_0^t K_v^{-1}\sigma(X_v^x,x)\cdot d\hat{B}_v$$
$$\gamma[X_t^x]=K_t[\int_0^tK_v^{-1}\sigma(X_v^x,x)\sigma^\ast(X_v^x,x)(K_v^{-1})^\ast dv]K_t^\ast$$
where $\sigma$ is the matrix whose columns are the $A_j$ $j=1,\ldots,d$ and $K$ the continuous invertible matrix valued process solution of 
$$K_t^x=I+\sum_{j=1}^d\int_0^t  \partial A^j(X_v^x,x)K_v^xdB_v^j+\int_0^t\partial B(X_v^x,x)K_v^xdv        .$$ where $\partial A^j$ and $\partial B$ are the Jacobian matrices with respect to the first argument.

We can write
$$
\Gamma[F]=\int_{[0,s]\times\mathbb{R}^m\times Y}\left(K^x_t\left[\int_0^t(K^x_v)^{-1}\sigma(X^x_v,x)\right.\right.\hspace{3cm}$$

$$\hspace{3cm}\sigma^*(X^x_v,x)(K^x_v)^{-1*}dv{\displaystyle]}(K^x_t)^*{\huge)}(y)M\odot\mu(dudxdy)$$

By the (EID) property, for $F$ to possess a density it suffises that the vector space $\mathcal{V}$ spanned by the column vectors of the matrices
$$\left(K_t^x(K_v^x)^{-1}\sigma(X_v^x,x)\right)(y)\qquad 0\leq v\leq t,\quad x\in \mathbb{R}^m,\quad y\in Y,$$
be $m$-dimensional  a.s. 

 If we restrict the study to the case where the diffusion coefficients do not depend on the first argument $A_j(X_u^x,x)=A_j(x)$, i.e. for the SDE
$$X_t^x=x+\sum_{j=1}^d A_j(x) B_t^j+\int_0^tB(X_v^x,x)dv$$
then, taking $v$ close to $t$, the space $\mathcal{V}$ contains the vectors
$$A_j(\Delta Z_u)\qquad j=1,\ldots,d\qquad u\in JT(Z)$$ where $JT(Z)$ denotes the jump times of $Z$ before $s$ and we have\\

\noindent{\bf Proposition 4}. {\it Let us suppose the L\'evy measure $\nu$ infinite. If the vectors $A_j(x)$ are such that for any infinite sequence $x_n\in\mathbb{R}^m$, $x_n\neq0$, tending to $0$, the vector space spanned by the vectors $$A_j(x_n),\;\; j=1,\ldots,d,\;\;n\in\mathbb{N} $$ is $m$-dimensional then the L\'evy process $(T_t(Z))_s$ has a density on $\mathbb{R}^m$.}\\

\noindent{\it Proof.} The result comes from the above condition by the fact that $Z$ has infinitely many jumps of size near zero
.\hfill$\Box$\\

As in part 2, the fact that $Z$ be a L\'evy process does not really matter. The method applies to the transform of the jumps of any process as soon as the perturbations are i.i.d and carry a Dirichlet form yielding (EID).

\end{document}